\documentclass{article}      
\usepackage{amssymb}
\usepackage{amsmath}
\usepackage{latexsym}
\usepackage[mathscr]{euscript}
\usepackage{graphicx}
\usepackage{amscd}
\usepackage{amsthm} 
\usepackage{fullpage} 
\usepackage{wrapfig}
\newtheorem{theorem}{Theorem}
\newtheorem{corollary}[theorem]{Corollary}
\newtheorem{lemma}[theorem]{Lemma}  
\newtheorem{proposition}[theorem]{Proposition}

\newtheorem{definition}{Definition}

\newcommand{\bz}{\mathbb{Z}}
\newcommand{\br}{\mathbb{R}}

\newcommand{\p}{\partial}

\newcommand{\ce}{\mathcal{E}}

\newcommand{\cg}{\mathcal{G}}
\newcommand{\cb}{\mathcal{B}}
\newcommand{\cp}{\mathcal{P}}

\newcommand{\cm}{\mathcal{M}}

\newcommand{\cf}{\mathcal{F}}

\newcommand{\cw}{\mathcal{W}}

\newcommand{\hk}{\hookrightarrow}

\newcommand{\med}{\medskip}
\newcommand{\la}{\longrightarrow}
\newcommand{\bfl}{\begin{flushleft}}
\newcommand{\efl}{\end{flushleft}}

\newcommand{\xr}{\xrightarrow}

\newcommand{\G}{\Gamma}

 \newcommand{\bcm}{\bar \cm}

 \newcommand{\cmlm}{\cm^{\ce}_\G(LM)}

  \usepackage{epsfig}

\DeclareMathOperator{\Crit}{Crit}

\DeclareMathOperator{\id}{id}

\DeclareMathOperator{\ind}{ind}

\def\R{\mathbb{R}}

\def\Z{\mathbb{Z}}

\begin{document}  

  \title{A Morse theoretic description of string topology  }
  \author{Ralph L. Cohen \thanks{The first author was partially
      supported by a grant from the NSF } \\ Department of Mathematics
    \\ Bldg. 380 \\Stanford University \\ Stanford, CA 94305, USA \and
    Matthias Schwarz \thanks{The second author was partially supported
      by a grant from the DFG}\\  Mathematics Institute \\ Universit\"at Leipzig\\ Postbox 10 09 20\\Leipzig D-04009  Germany
  }
\date{\today}
\maketitle  

\begin{centerline} {\sl This paper is dedicated to Yasha Eliashberg on the occasion of his 60th birthday.}  \end{centerline}

\rm
 \begin{abstract}  Let $M$ be a closed, oriented, $n$-dimensional manifold. In this paper  we give a Morse theoretic description of the string topology operations introduced by Chas and Sullivan, and extended by the first author, Jones, Godin, and others.  We do this by studying maps  from surfaces with cylindrical ends to $M$, such that on the cylinders, they satisfy the gradient flow equation of a Morse function on the loop space, $LM$.   We then give Morse theoretic descriptions of  related constructions, such as the Thom and Euler classes of a vector bundle, as well as the shriek, or umkehr homomorphism.     \end{abstract}

  \tableofcontents

 \section*{Introduction}

String topology operations were first  defined by Chas and Sullivan in \cite{CS}. Their basic loop product is an algebra structure on the homology of the free loop space of a  closed, oriented manifold, $H_*(LM)$,  which comes from studying maps from a ``pair of pants" surface to $M$.   Since then generalizations and applications of these operations have been widely studied.  In particular  the work of V. Godin  \cite{godin} describes operations based on families of Riemann surfaces, varying in moduli space. 

More specifically, 
let $\cm_{g,p+q}$ be the space of oriented, connected surfaces embedded in $\br^\infty$ having  genus $g$ and  $p+q$ parameterized boundary components. We think of these surfaces as  cobordisms between $p$ parameterized circles,  thought of as ``incoming", 
 and $q$ parameterized circles,  thought of as ``outgoing".  This space is homotopy equivalent to the moduli space of bordered Riemann surfaces, with marked points on each boundary component,  and is a model for the classifying space, $BDiff^+ (\Sigma_{g, p+q}; \p \Sigma)$, where $Diff^+ (\Sigma_{g, p+q}; \p \Sigma)$ is the group of orientation preserving diffeomorphisms of a surface $\Sigma$ that are fixed pointwise on the boundary. 
 
  For a closed, oriented $n$-manifold $M$, let $\cm_{g,p+q}(M)$ denote the space of pairs,
  $$
\cm_{g,p+q}(M) = \{(\Sigma, f): \, \Sigma \in \cm_{g,p+q} \, \text{and} \, f : \Sigma \to M \, \text{is a smooth map} \}.
$$

$\cm_{g,p+q}(M)$ is a model of the homotopy orbit space,
$$
\cm_{g,p+q}(M) \simeq EDiff^+ (\Sigma_{g, p+q}; \p \Sigma) \times_{Diff} C^\infty (\Sigma, M)
$$
where the subscript $Diff$ refers to taking the orbit space by the diagonal $Diff^+ (\Sigma_{g, p+q}; \p \Sigma) $-action.

     \begin{figure}[ht]
  \centering
 \includegraphics[height=7cm]{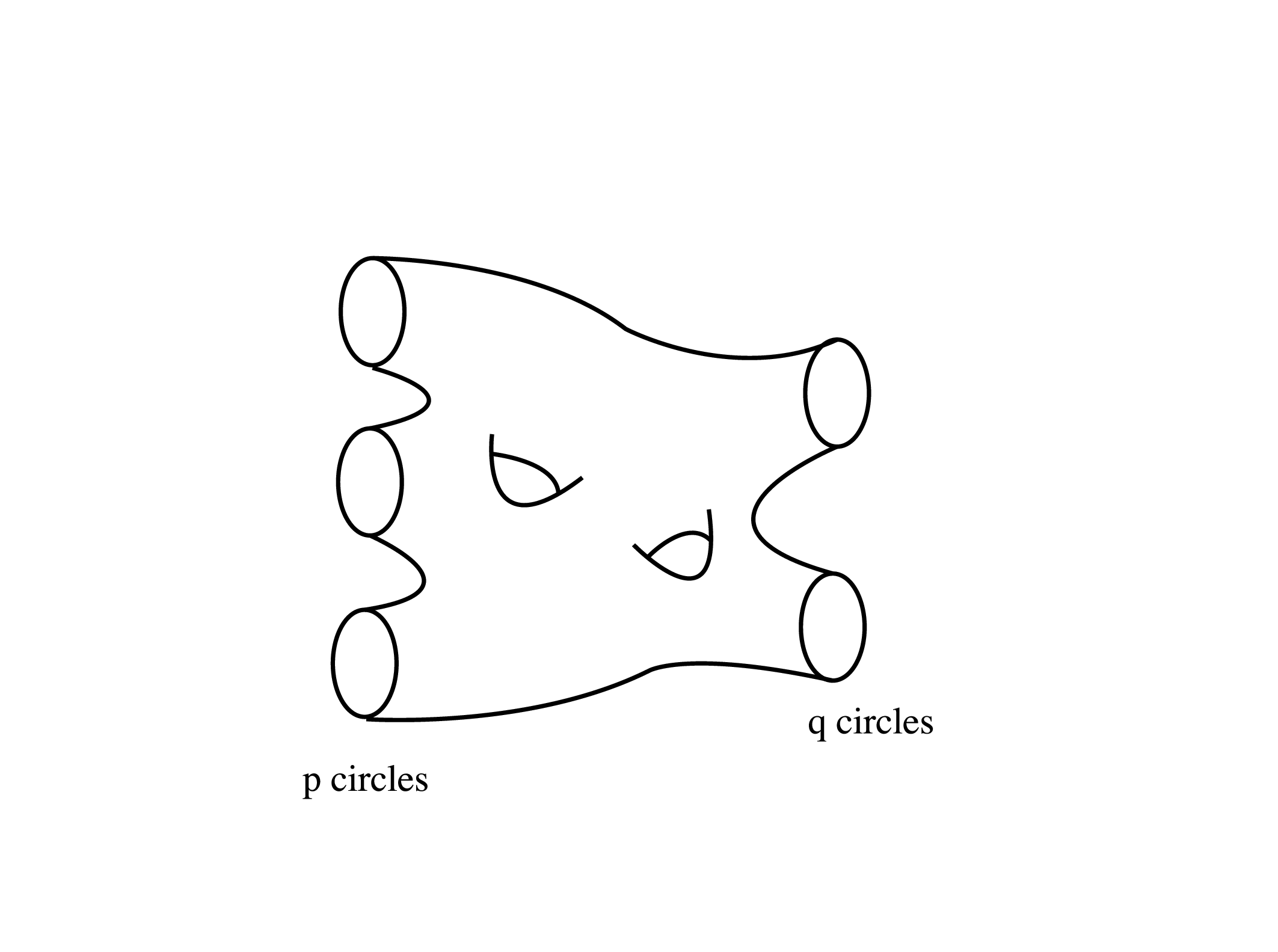}
\caption{The  surface $\Sigma$}
 \end{figure}

Using the restriction maps to the incoming and outgoing parameterized boundary components as well as the projection map $\cm_{g,p+q}(M) \to \cm_{g,p+q}$, one has a diagram,   

\begin{equation}\label{restrict}
(LM)^q \xleftarrow{\rho_{out}} \cm_{g,p+q}(M) \xr{\rho_{in}} \cm_{g,p+q} \times (LM)^p.
\end{equation}
 In \cite{godin},  this correspondence diagram lead to   higher string topology operations in a (generalized) homology theory $h_*$, by first constructing 
an  ``umkehr map"
$$
(\rho_{in})_! : h_*(\cm_{g,p+q} \times (LM)^p) \to h_{* +\chi (F)\cdot d}(\cm_{g,p+q}(M))
$$ where $\chi (F)$ is the Euler characteristic of the surface $\Sigma$ and $d = dim (M)$.   For this one needs that the generalized homology theory supports an orientation of $M$.    The higher string topology operations of \cite{godin} were then defined via the composition
\begin{equation}\label{operate}
\mu_{g, p+q} : h_*(\cm_{g,p+q} \times (LM)^p) \xr{(\rho_{in})_!} h_{* +\chi (\Sigma)\cdot d}(\cm_{g,p+q}(M)) \xr{(\rho_{out})_*} h_{*+\chi (\Sigma)\cdot d}((LM)^q).
\end{equation}

In the case of an operation based on a fixed surface (i.e a point in $\cm_{g,p+q}$), and in the  case of families of genus zero surfaces when $q = 1$,  the
 umkehr map  has been given several descriptions.  In \cite{CS} it was constructed on the chain
level, and in \cite{CJ} and \cite{CG} it was constructed as a generalized Pontrjagin-Thom construction on the homotopy theoretic level.   Pontrjagin-Thom constructions are also the basis of the umkehr map defined by Godin in \cite{godin} for surfaces varying in higher genus moduli space.  

\med
 In this paper we will  show how  string topology operations can be constructed using Morse theory on the   loop space $LM$.     In section one we show how to construct the umkehr map, and therefore the resulting string topology operations on the level of the Morse chain complex of appropriate energy functions on $LM$.  We will prove that the operations defined this way are equal  to the original string topology operations on the level of homology.   We then  indicate how this construction can be generalized to families, using the work of Godin \cite{godin}.        In section 2  we will show that  under the appropriate transversality conditions, operations obtained 
   by explicitly counting the ``gradient flow surfaces"   whose boundaries
 lie in appropriate stable and unstable manifolds of critical points,  define the same string topology operations.  We discuss these transversality conditions in some detail.  The operations constructed this way were defined and studied by the second author with Abbondandolo in \cite{abboschwarz2}.   In that paper the authors described an  isomorphism  of rings  between the Floer   homology of the cotangent bundle, $HF_*(T^*M)$ with the ``pair of pants" product, to this
 Morse loop product in  $H_*(LM)$.  
 The arguments in this paper verify not only in full detail that this product agrees with the string topology product as constructed by Chas and Sullivan \cite{CS}, but also give this equivalence for more general chord diagrams and operations.  This verification uses related Morse theoretic constructions of the Thom and Euler classes of any oriented vector bundle, as well as  the  ``umkehr" map, which is done  in section 3.   These constructions may be of independent interest. 
 
  \section{Fat graphs and the space of gradient surfaces in a manifold}
 We begin by describing the types of Morse functions on the loop space that we will consider.  We refer the reader to \cite{abboschwarz} for more details.
 
Endow $M$ with a Riemannian metric. Consider  a smooth Lagrangian
$$
L : \br/\bz \times TM \to \br
$$
that satisfies the following   convexity property, bounds on its second derivatives, as well as   nondegeneracy properties.    

\bfl 
\bf (L1) \rm There exists $\ell_0 > 0 $ such that

\efl
$$
\nabla_{v,v}L(t, (q,v)) \geq \ell_0I
$$
for every $(t, (q,v)) \in  \br/\bz \times TM$.  (Here $q \in M$ and $v \in T_qM$.)

\bfl 
\bf(L2) \rm There exists $\ell_1 \geq 0$ such that

\efl
$$
|\nabla_{v,v}L(t, (q,v))| \leq \ell_1 \quad |\nabla_{q,v}L(t, (q,v))| \leq \ell_1(1 + |v|), \quad |\nabla_{q,q}L(t, (q,v))| \leq \ell_1(1 + |v|^2)
$$
for every $(t, (q,v)) \in \br/\bz \times TM$. 

\med
We explain a bit of this notation.  The Riemannian metric on $M$ induces a splitting of the tangent bundle $T(TM)$ into a vertical and horizontal part, via the Levi-Civita connection.  Then $ \nabla_{v,v}$, $ \nabla_{q,v}$, and $ \nabla_{q,q}$ denote the components of the Hessian in this splitting.   

With such a Lagrangian one can define an energy function,
 
\begin{align}
\ce : LM &\la \br   \notag \\
\ce (\gamma) &= \int_0^1 L(t, \gamma (t), \frac{d \gamma}{ds}(t)) dt  \notag
\end{align}
  which is $C^2$ on the space $LM$, which we take to be those loops of Sobolev class $W^{1,2}$. 

However, if we want the energy function to be smooth on this Hilbert manifold of loops, we have to assume stronger conditions on the Lagrangian $L$, namely similar bounds on all partial derivatives of $L$. For example, for any Riemannian metric $g$ on $M$ and a time-dependent potential $V(t,q)$ on $M$ we can take $L(t,(q,v))=\frac{1}{2}|v|_g^2 + V(t,q)$, and $\ce$ will be smooth.  In Section 2, it will be necessary to assume such smoothness
 of $\ce$ for reasons of transversality.  Therefore, we will from now on only consider such Lagrangians with a fibrewise quadratic kinetic term. 
 However, the final result on the equivalence of the Morse-theoretic construction of the string topology operation and the homotopy-theoretic definition is valid for the more general Lagrangians satisfying (L1) and (L2) due to the canonical continuation homomorphism between the Morse homologies for different Lagrangians.

   We then also assume
  
  \bfl
  \bf (L0) \rm The critical points of $\ce$ are all nondegenerate.
  
  \efl
  We denote the set of critical points of $\ce$ by $\cp (L)$.

In this context, the energy functional,  $\ce : LM \la \br$ 
 is a Morse function that is bounded below with critical points of finite Morse indices,  and it satisfies the Palais-Smale condition.  Again, we refer the reader to \cite{abboschwarz} for details.

When the Lagrangian $L$ satisfies these assumptions, standard Morse theory applies, and one can construct a  space $LM_\ce$  which is defined to be the union of the unstable manifolds of $\ce$, and is topologized as a subspace of the loop space, $LM$.  $LM_{\ce}$    has one cell for each critical point in $\cp(L)$, and  the inclusion $LM_{\ce} \hk LM$ is a homotopy equivalence.  The cellular chain complex of $LM_\ce$ is the Morse complex,
\begin{equation}\label{complex}
\la \cdots \xr{\p_{p+1}} C_p^\ce (LM) \xr{\p_p} C^\ce_{p-1}(LM) \to \cdots
\end{equation}
where $C_p^\ce(LM)$ is the free abelian group generated by those $a \in \cp(L)$ of Morse index $p$, and 
$$
\p_p([a]) = \sum_{\substack{b\in \cp(L)\\ ind (b) = p-1}} \#\cm(a,b) [b]
$$
where 
 $\cm(a,b)$ is the space of gradient flow lines connecting $a$ to $b$, which
 is a compact, zero dimensional, oriented manifold in this setting.  The number $ \#\cm(a,b)$ refers to the oriented count of the points in this moduli space which requires a choice of orientations for each cell
 We note that the stable attaching maps of the cells of the   $LM_\ce$ can be described by the framed bordism types of the higher dimensional compact  spaces of  piecewise flows, $\bcm (a,b)$ \cite{cjs}\cite{cfloer}.
 
 \med
 The homotopy theoretic  string topology operations were defined using ``fat graph" models for surfaces \cite{CG} \cite{godin}.  We will likewise use these graphs to define our Morse theoretic operations.

 \med
 
 We recall the definition ( see \cite{penner}, \cite{strebel}).
\begin{definition}   A fat graph is a finite graph with the following properties:
\begin{enumerate} 
\item Each vertex is at least trivalent 
\item Each vertex comes equipped with a cyclic order of the half edges emanating from it.
\end{enumerate}
\end{definition}

The cyclic order of the half edges is quite important in this structure.
It allows for the graph to be ``thickened" to a surface with boundary.  As a way of describing this thickening,  recall that the cyclic orderings of the half edges at each vertex   define a  partition of the set   of oriented edges, that identify boundary components of the thickened surface.  More explicitly,    let $E(\G)$ be the set of edges, and let $\tilde E(\G)$ be the set of oriented edges.  $\tilde E(\G)$ is a $2$-fold cover of $E(\G)$.  It has an involution $e \to \bar e$
which represents changing the orientation.  The partition of $\tilde E(\G)$ is best illustrated by the following example.   

\begin{figure}[ht]
 \centering
  \includegraphics[height=7cm]{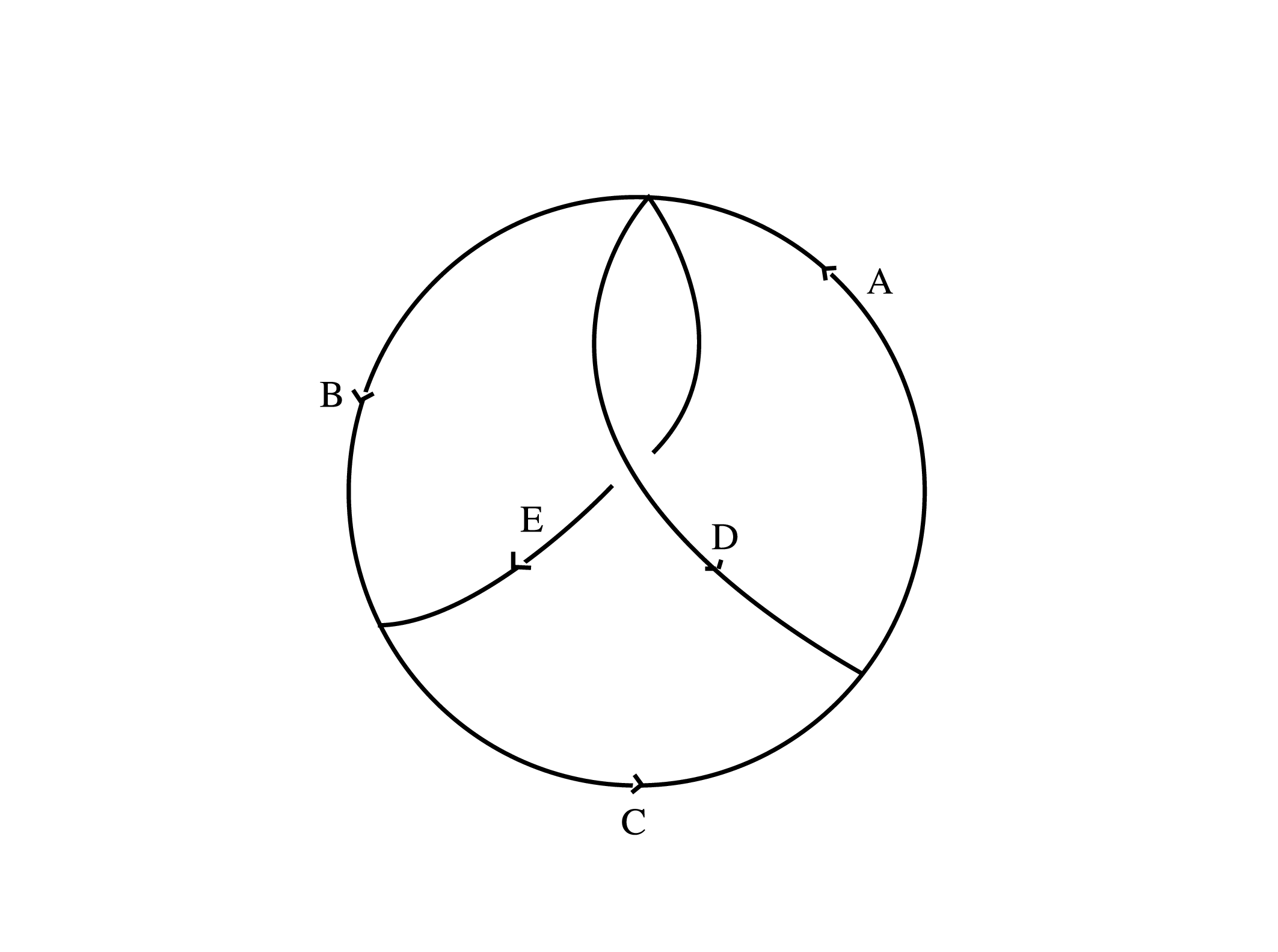}
 \caption{The  fat graph  $\G$}
   \label{fig:figtwo}
\end{figure}
 
 In this picture the cyclic   orderings  at the vertices are determined by the counterclockwise orientation of the plane. 
 To obtain the partition,  notice that an oriented edge has   well defined source and target vertices. Start with an oriented edge, and follow it to its target vertex.  The next edge in the partition is the next oriented edge in the cyclic ordering at that vertex.  Continue in this way until one is back at the original oriented edge.  This will be the first cycle in the partition.  Then continue with this process
 until one exhausts all the oriented edges.  The resulting  cycles in the partition will be called  \sl ``boundary cycles" \rm as they reflect the boundary circles of the thickened surface. In the  case of $\G_2$ illustrated in figure 2, the partition into boundary cycles are given by:
  
  $$
  \text{Boundary cycles of $\G_2$:} \quad (A,B,C) \, (\bar A, \bar D, E, \bar B, D, \bar C, \bar E ).
  $$
  
 So one can  read off 
 the number of boundary components in the thickened surface of a fat graph.  Furthermore the graph and the surface have the same homotopy type, so one can compute the Euler characteristic of the surface directly from the graph.  Then using the formula
 $\chi (F) = 2 - 2g -n$, where $n$ is the number of boundary components, we can solve for the genus directly in terms of the graph.  The space of \sl metric fat graphs \rm (i.e fat graphs with lengths assigned to each edge)  of topological type $(g,n)$ gives a model for the homotopy type of the 
 moduli space of Riemann surfaces of genus $g$ and $n$-marked points \cite{penner}, \cite{strebel}.  
 
  Notice that the boundary cycles of a metric fat graph $\G$ nearly determines a parameterization of the boundary of the thickened surface.  For example, the boundary cycle $(A, B, C)$ of the graph $\G $ above  can be represented by a map $S^1 \to \G_2$ where the circle is of circumference equal to the sum of the lengths of sides $A$, $B$, and $C$.
  The ambiguity of the parameterization is the choice of where to send the basepoint $1 \in S^1$.  By choosing a marked point in each boundary cycle,  Godin   \cite{godin1} \cite{godin}   used ``marked" fat graphs to give   models of the homotopy type of  moduli spaces of bordered Riemann surfaces.    All the fat graphs we will work with will be assumed to be marked in this way.   Also, as part of our data in a marking of a metric fat graph, we assume that the $n$ boundary cycles  are partitioned into $p$ incoming, and $q = n-p$ outgoing cycles.
  
  Let $\G$ be a metric marked fat graph.   In particular    this means that the boundary cycles of $\G$ are partitioned into $p$ incoming and $q$ outgoing cycles, and there are parameterizations determined by the markings,
\begin{equation}\label{alphas}
  \alpha^- : \coprod_p S^1 \la \G, \quad \alpha^+ : \coprod_q S^1 \la \G.
\end{equation}
  
  By taking the circles to have circumference equal to the sum of the lengths of the edges making up the boundary cycle it parameterizes, each component of $\alpha^+$ and $\alpha^-$ is a local isometry.   
  
  Define the surface $\Sigma_\G$ to be the mapping cylinder of these parameterizations,
  \begin{equation}\label{sigmag}
  \Sigma_\G = \left( \coprod_p S^1 \times (-\infty, 0]\right)   \sqcup \left(\coprod_q S^1 \times [0, +\infty)\right)  \bigcup  \G / \sim
  \end{equation}
  where $(t,0) \in S^1 \times (-\infty, 0] \, \sim \alpha^- (t) \in \G$, and $(t,0) \in S^1 \times [0, +\infty)  \, \sim \alpha^+ (t) \in \G$.
  
  \med
  
\begin{figure}[ht]
  \centering
  \includegraphics[height=8cm]{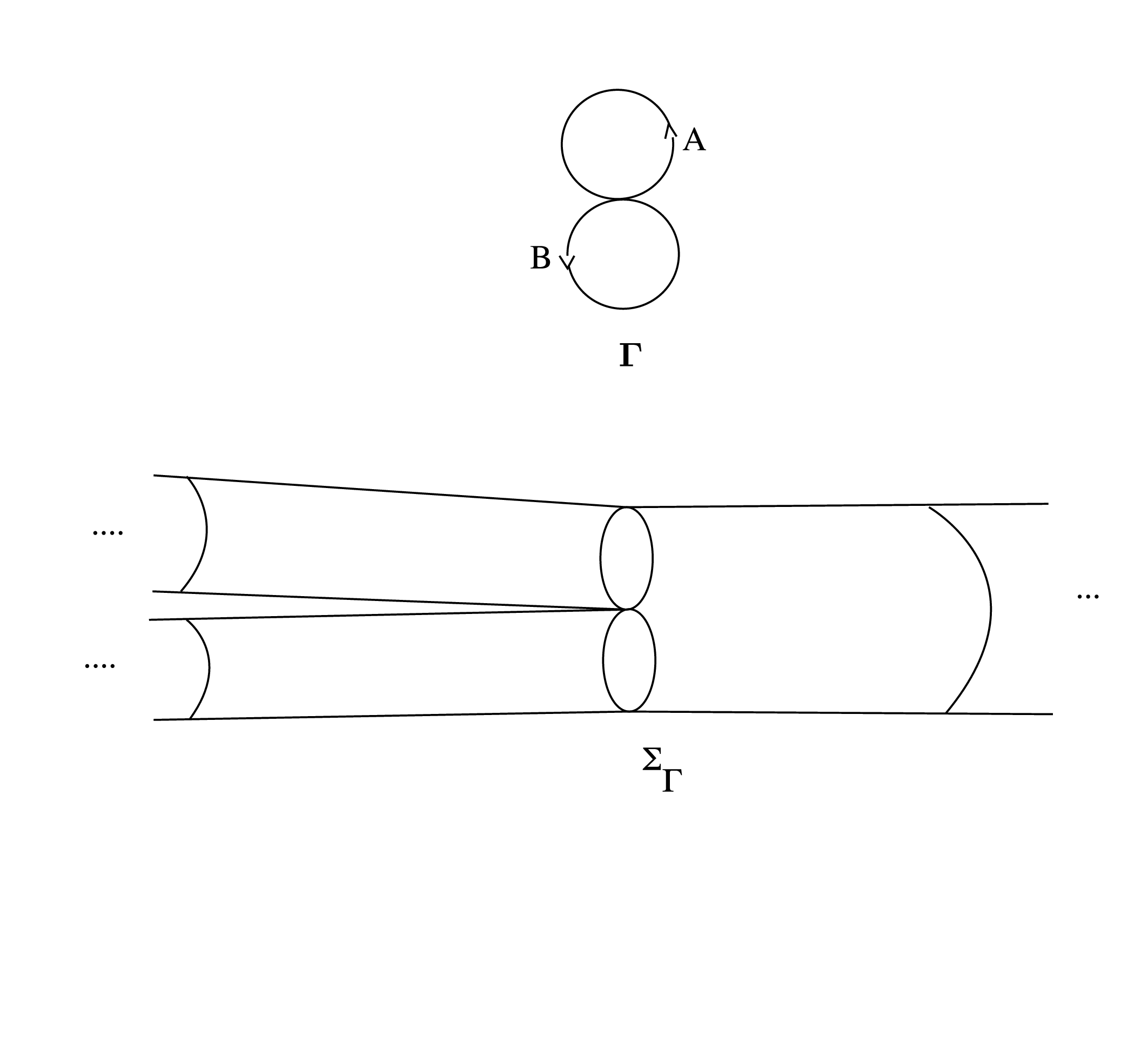}
 \caption{$\Sigma_\G$ }
   \label{fig:figthree}
\end{figure}

  Notice that the graph in  figure 3  is a  fat graph representing a surface of genus $g=0$ and $3$ boundary components.  This graph has two edges, say $A$ and $B$, and has boundary cycles
   $(A), (B), (\bar A, \bar B)$.  If we let $(A)$ and $(B)$ be the incoming cycles and $(\bar A, \bar B)$ the outgoing cycle, then figure \ref{fig:figthree}  is a picture of the surface $\Sigma_\G$, for $\G$ equal to the figure 8.

Notice that a map    $\phi : \Sigma_\G \to M$ is  a collection of $p$ maps from half cylinders,  $\phi_i : (-\infty, 0] \times  S^1 \to  M$ and $q$- half cylinders,
$\phi_j : [0, +\infty) \times S^1  \to  M$, that have an intersection property at $t=0$ determined by the combinatorics of the fat  graph $\G$.
In the definition of $\phi_i : (-\infty, 0) \times S^1 \to M$, the circle factor is rescaled in a canonical way so as to have radius one.

We now define a ``gradient flow surface" to be a map $\phi : \Sigma_\G \to M$ that when restricted to each half cylinder satisfies  a gradient flow equation.  Now, similar to what was done in constructing cohomology operations on closed manifolds using Morse theory (see \cite{betzcohen}, \cite{fukaya}, \cite{morsefield}),  we will need to allow our gradient flow equations to be perturbed on each cylinder.  More specifically, define a ``Lagrangian labeling"  of a marked fat graph $\G$ to be a labeling $\ce (\G)$  of each of the  boundary cycles of $\G$ by a Lagrangian, $L_i$, and therefore by an energy functional $\ce_i : LM \to \br$.  We write
$\ce (\G) = (\ce_1, \cdots \ce_{p+q})$   where $\ce_i$ is the energy functional labeling the $i^{th}$ boundary cycle. 

\med
\begin{definition}\label{cmlm}  Let $\G$ be a marked fat graph with Lagrangian labeling $\ce = \ce (\G)$.    Define the moduli space of ``gradient flow surfaces", $\cm^{\ce}_\G(LM)$ to be the space of maps
$$
  \phi:  \Sigma_{\G} \to M
$$
  that  are  smooth in the interiors of the  cylinders,  and that the restrictions to the  incoming  cylinders define maps  $\phi_i :  (-\infty, 0) \times   S^1 \to  M$, $i = 1, \cdots , p$   satisfying
  the gradient flow equation
\begin{equation}\label{gradflow}
  \frac{d\phi_i (t,s)}{dt}  + \nabla \ce_i  = 0
\end{equation}
%
  such that $\lim_{t \to -\infty} \phi_i(t,\cdot) : S^1 \to M$  converges uniformly to  a critical point in $\cp (L_i)$.
  Similarly on outgoing cylinders $\phi$ defines maps $\phi_j : (0,+\infty) \times S^1 \to M$, $j = 1, \cdots,  q$
  which satisfy the gradient flow equation  $  \frac{d\phi_j (t,s)}{dt}  + \nabla \ce_{p+j}  = 0$ and $\lim_{t\to +\infty} \phi_j(t, \cdot) : S^1 \to M$ converges uniformly to a critical point in $ \in \cp (L_{p+j})$.    
  
  The space $\cm^{\ce}_\G(LM)$ is topologized as a subspace of the space of continuous maps $\Sigma_{\G} \to M$  in the compact-open topology.
  \end{definition}
  
  \med
  The spaces $\cmlm$ will be essential in our definition of the Morse theoretic string    topology.   For example, we now
  describe a correspondence diagram analogous to (\ref{restrict}).
  
   Let $\phi \in \cmlm$.  For $i = 1, \cdots , p$, let $\phi_{i,-1} : S^1 \to M$ be the restriction of $\phi_i : S^1 \times (-\infty, 0] \to M$ to $S^1 \times \{-1\}$.  Notice that by definition, each $\phi_{i,-1}$ lies in an unstable manifold of some critical point in $\cp (\ce_i)$.  Therefore $\phi_{i,-1} \in LM_{\ce_i}$.  
   
   Similarly, for $j = 1, \cdots, q$,
  let $\phi_{j,1} : S^1 \to M$ be the restriction of $\phi_j$ to $S^1 \times \{1\}$.  These restrictions define  the following maps. 
  
  \begin{equation}\label{restrict2}
 \prod_{j=1}^q LM_{\ce_{p+j}}\xleftarrow{\rho_{out}}  \cmlm \xr{\rho_{in}}  \prod_{i=1}^p LM_{\ce_i}.
  \end{equation}
  Our goal is to construct an umkehr map on the level of chains,
  $$
  (\rho_{in})_! :  \bigotimes_{i=1}^p C_*^{\ce_i}(LM)   \to C_*(\cmlm)
  $$
  so that our  string topology operation on the level of Morse homology will be  induced by the composition on the level of Morse chains,
  $$
  \mu_\G:   (\rho_{out})_* \circ (\rho_{in})_! :  \bigotimes_{i=1}^p C_*^{\ce_i}(LM)    \la  \bigotimes_{j=1}^q C_*^{\ce_{p+j}}(LM).
  $$
 \med
 To do this, we will find it convenient, as was the case in 
  \cite{CG},   to consider a particular type of marked fat graph, known as a ``Sullivan chord diagram".
  
  \med 

\begin{definition}   A  ``Sullivan chord diagram"  of type $(g; p,q)$ is a  fat graph  representing a surface of genus $g$ with
$p+q$ boundary components, that consists of  a disjoint union of 
$p$  disjoint closed circles   together with the disjoint union of connected trees
whose endpoints lie on the circles.    The cyclic orderings  of the edges at the vertices must be such that each of the $p$ disjoint
circles is a boundary cycle.   These $p$ circles are referred to as the incoming boundary cycles, and the other $q$ boundary cycles
are referred to as the outgoing boundary cycles. 
\end{definition}

The ordering at the vertices in the diagrams that follow
are indicated by the  counterclockwise cyclic ordering of the plane
Also in a Sullivan chord diagram, the vertices and edges that lie on one of the $p$ disjoint circles will
be referred to as circular vertices and circular edges respectively.  The others will be referred
to as ghost vertices and edges.   From now on we refer to a ``Sullivan chord diagram" as simply a ``chord diagram".  
 
\begin{figure}
\begin{center}
 \includegraphics[height=5cm]{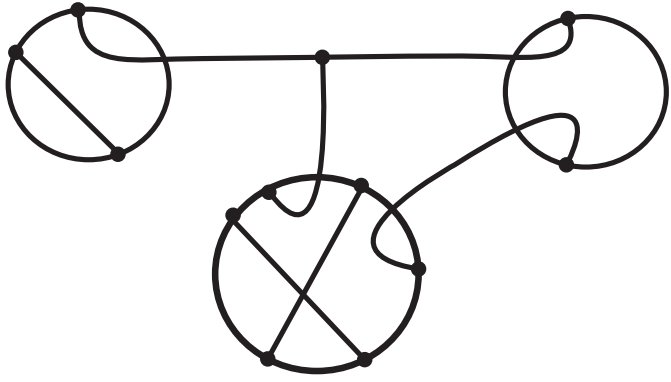}
 \caption{Sullivan chord diagram of type (1;3,3)}
\end{center}
\end{figure}

Let $\G$ be a marked chord diagram
The marking defines a parameterization
of the incoming and outgoing boundary circles, and hence if $\phi : \G \to M$, we can identify the restriction to these boundary circles with loops, $\phi_{i,0} : S^1 \to M$.

We now go about studying the topology of $\cmlm$.  The following is our main result.

\med

\begin{theorem}\label{main}  Let $\G$ be a marked chord diagram. 
The natural map from the space of gradient flow surfaces to the continuous mapping space,
$$\cmlm \to Map (\Sigma_\G, M)$$
is a homotopy equivalence.
\end{theorem}

\med
\begin{proof}
Let $\G$ be a chord diagram, and let  $v(\G)$ be the collection of circular vertices (i.e vertices that lie on the incoming
boundary circles).  There is a natural evaluation map
$$
ev_\G :  \prod_{i=1}^p LM_{\ce_i} \to M^{v(\G)}
$$  that evaluates the $i^{th}$ loop on the vertices lying on the $i^{th}$ boundary circle of $\Gamma$.  
Put an equivalence relation on the set of circular  vertices $v(\G)$ by saying that two vertices $v_1$ and $v_2$ are equivalent if
there is a ghost subtree of $\G$ that contains both $v_1$ and $v_2$.  Let $\sigma(\G) = v(\G)/\sim$ be the set of equivalence classes of these circular vertices.  The projection map 
$\pi : v(\G) \to \sigma (\G)$ defines a diagonal embedding  

\begin{equation}\label{diag}
\Delta_\G: M^{\sigma (\G)} \hk M^{v(\G)}.
\end{equation}

\begin{lemma}\label{lemma}
$\cmlm$ is homotopy equivalent to the homotopy pullback of the map $
ev_\G :  \prod_{i=1}^p LM_{\ce_i} \to M^{v(\G)}$ along the diagonal embedding $\Delta_\G: M^{\sigma (\G)} \hk M^{v(\G)}$.
\end{lemma}
\begin{proof}    
We now describe a locally trivial fiber bundle,
$$
{\tilde{ev_\G}} :   (LM_{\ce})^p_{\G} \to M^{v(\G)}
$$
of the same homotopy type as the evaluation map, $ev_\G :  \prod_{i=1}^p LM_{\ce_i} \to M^{v(\G)}$.  

We define $(LM_\ce)^p_\G$ to be the space of ``hairy loops" defined by the graph $\G$.  Namely,
let $C_1, \cdots C_p$ be the $p$ -incoming circles of  chord diagram $\G$.  Let $v_{i,1},  \cdots v_{i, n_i} \subset C_i$
be the set of circular vertices lying on $C_i$.  We define ``hairy incoming circles" by attaching intervals at these vertices:  Let
$C^h_i = C_i \cup \bigcup_{n_i}[0,1]$ where the $j^{th}$ interval is attached at $t = 0$ to the $j^{th}$ vertex $v_{i,j} \in C_i$. 
We now define the space of hairy loops as follows.  We let $(LM_\ce)^p_\G =  \{ \theta \in Map (\cup_{i=1}^p C^h_i,  M) \, :  \, \theta_{|_{C_i}} : C_i \cong S^1 \to M  \, \text{lies in} \, LM_{\ce_i} \}.
  $ We have an inclusion $\iota :   \prod_{i=1}^p LM_{\ce_i} \hk (LM_\ce)^p_\G$ which are maps that are  defined to be constant on the intervals (``hairs").   Clearly this map is a homotopy equivalence.

 \begin{figure}
\begin{center}
 \includegraphics[height=8cm]{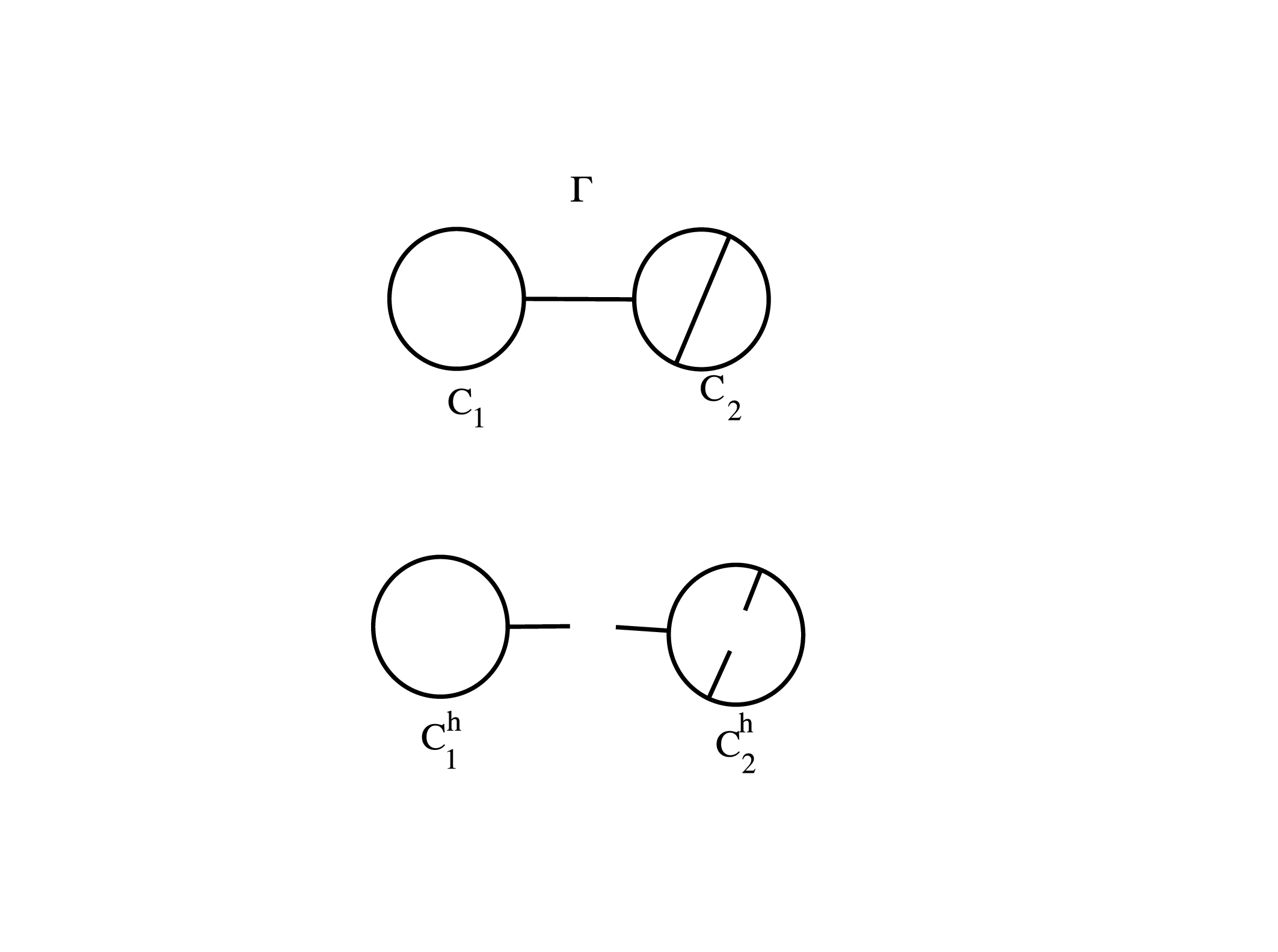}
 \caption{The hairy incoming circles of a chord diagram $\G$}
\end{center}
\end{figure}

Notice that another way to describe the space hairy loops is as follows:
$$
(LM_\ce)^p_\G = \{(\gamma, \alpha) \in   \prod_{i=1}^p LM_{\ce_i} \times (M^{v(\G)})^I: \, ev_\G (\gamma) = \alpha (0) \}
$$
where $X^I$ denotes the space of paths $\alpha : I = [0,1] \to X$.  This implies that there is a Serre fibration,
\begin{align}\label{fibration}
\tilde{ev}_\G : (LM_\ce)^p_\G  &\la M^{v(\G)} \\
(\gamma, \alpha) &\la \alpha (1) \notag
\end{align}
Indeed this fibration has the structure of a locally trivial fiber bundle.  (See \cite{kling} for descriptions of local trivializations.)

Let $P_\G$ be the pullback (restriction) of the bundle $\tilde{ev}_\G : (LM_\ce)^p_\G  \la M^{v(\G)}$ to the image of the embedding, $\Delta_\G (M^{\sigma (\G)}) \subset M^{v(\G)}$

\begin{equation}\label{pullback}
\begin{CD}
P_\G  @>\hk >>   (LM_\ce)^p_\G   \\
@V\tilde{ev} VV     @VV\tilde{ev}V \\
M^{\sigma (\G)}   @>\hk >\Delta_\G > M^{v(\G)}
\end{CD}
\end{equation}

Now by definition, the $P_{\G}$ is defined to be the space of pairs $((\gamma, \alpha) \in  (LM_\ce)^p_\G$ such that  $ \alpha (1) \in \Delta_\G (M^{\sigma (\G)})$.   This space can be described alternatively as follows.

Let $\tilde \G$ be the graph constructed from the union of the hairy incoming circles,

$$
\tilde \G = \bigcup_{i=1}^p  C^h_i /\sim
$$
where we make the following identifications:   We identify the endpoint of the ``hair"  (i.e  $t=1$ in the interval) emanating from vertex $v_1$ with the endpoint of the hair emanating from vertex $v_2$ if and only if $v_1$ and $v_2$ are vertices of the same ghost subtree in $\G$.  That is, they are identified if and only if these vertices are in the same equivalence class with respect to the relation
defined in (\ref{diag}) above.  We initially put a vertex at each of these identification points, but then remove those new vertices  that are only bivalent.  

\begin{figure}
\begin{center}
 \includegraphics[height=5cm]{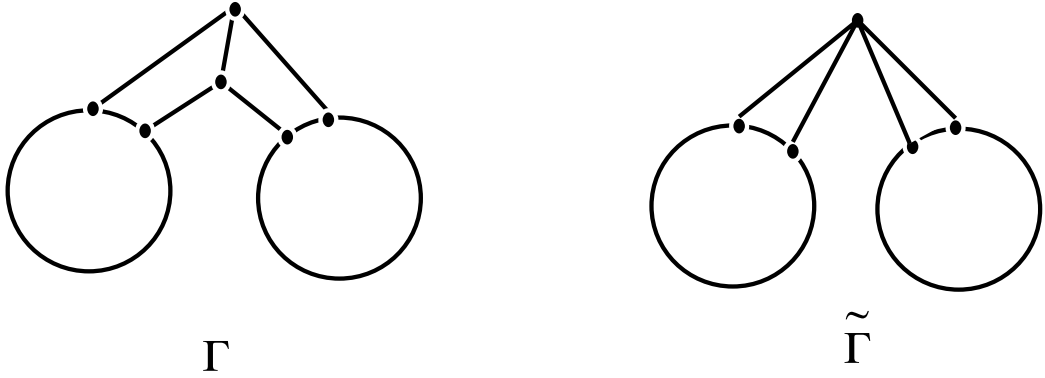}
 \caption{The chord diagrams $\G$ and $\tilde \G$}
\end{center}
\end{figure}

Notice then that
  $P_\G$ consists of maps $\theta : \tilde \G \to M$  whose restriction
to the $i^{th}$ incoming circles, $C_i$ lies in in $LM_{\ce_i}$,  $i = 1, \cdots , p$.  Notice furthermore, that there is a natural map of graphs 
$$
p : \G \to \tilde \G
$$
that is defined by collapsing various ghost trees in $\G$ to the new chord vertices of $\tilde \G$ defined above.    In particular, $p : \G \to \tilde \G$ is a homotopy equivalence.
Moreover, the cyclic ordering of the half edges emanating at the  vertices in $\G$ defines a cyclic ordering of the half edges emanating at the vertices of $\tilde \G$, and so
$\tilde \G$ has the structure of a fat graph, and indeed a Sullivan chord diagram
with the same marked incoming circles as $\G$.  

Let $ Map_\ce (\G, M)$ denote the space of maps $\beta : \G \to M$ whose restriction to the $i^{th}$ incoming circle lies in $LM_{\ce_i}$.  The map $p : \G \to \tilde \G$ defines a homotopy equivalence $P_\G 
  \simeq Map_\ce (\G, M)$.  But this latter space
is homeomorphic to $\cmlm$.  This can be seen as follows.  Since   any
map $\theta : \G \to M$ whose restrictions to the incoming circles are
loops in $\theta_i \in LM_{\ce_i}$, each of these loops extends in a
unique way  to a map  of a  half cylinder  $\bar \theta_i : (-\infty,
0] \times S^1 \to M$    satisfying the gradient flow equation.
(\ref{gradflow})
 Moreover, since $\ce$ satisfies the Palais-Smale criterion and is bounded below, the restrictions of $\theta$ to the outgoing
boundary cycles, $\theta_j : S^1 \to M$ also extend to a gradient flow cylinder, $\bar \theta_j : [0, +\infty) \times S^1 \to M$.  These gradient flow cylinders patch together to give an element in $\cmlm$.   Thus  $P_\G \simeq \cmlm$, which proves the lemma.
\end{proof}

We now complete the proof of Theorem \ref{main}. Since $\ce : LM \to \br$ is Palais-Smale and bounded below, the inclusion $LM_\ce \hk LM$ is a homotopy equivalence. Thus the homotopy pullback of $ev_\G :  \prod_{i=1}^p LM_{\ce_i} \to M^{v(\G)}$ along the diagonal embedding $\Delta_\G: M^{\sigma (\G)} \hk M^{v(\G)}$,  is homotopy equivalent to the pullback of the fibration $ev_\G : LM^p \to M^{v(\G)}$ along the diagonal $\Delta_\G$.  Now as described in \cite{CG}, this pullback is the mapping space, $Map (r(\G), M)$, where $r(\G)$ is the ``reduced"  chord diagram obtained from $\G$ by collapsing each ghost edge
to a point.  By the lemma, we then have a homotopy equivalence, $\cmlm \simeq Map (r(\G), M)$.   Now since the collapse map $\G \to r(\G)$ is a homotopy equivalence, we have an equivalence of mapping spaces, $Map (r(\G), M) \simeq Map (\G, M)$.  But this last mapping space is homotopy equivalent to $Map (\Sigma_\G, M)$, since the surface $\Sigma_\G$ retracts onto the graph $\G$.
This completes the proof of the theorem.
\end{proof}
Notice that this argument yields a commutative diagram
$$
\begin{CD}
P_\G  @<\simeq <<  Map_\ce (r(\G), M)@>\hk >>  \prod_{i=1}^p LM_{\ce_i} \\
&& @VVV   @VVV \\
&& Map (r(\G), M) @>\hk >> (LM)^p
\end{CD}
$$
where the two vertical maps are homotopy equivalences, as well as the upper left horizontal map.    In particular $Map_\ce (r(\G), M) \hk  \prod_{i=1}^p LM_{\ce_i}$ and $Map (r(\G), M)  \hk  (LM)^p$ are both topological embeddings with open neighborhoods  given by the inverse image of a tubular neighborhood $\eta (\Delta_\G)$ of the embedding $\Delta_\G : M^{\sigma (\G)}    \hk > M^{v(\G)}$ of compact manifolds. Even though $Map_\ce (r(\G), M)$ is not smooth, we think of these neighborhoods as ``tubular neighborhoods",  since they are   
  homeomorphic to the total spaces of the pullbacks of the normal  bundle,   $\nu (\Delta_\G) \to M^{\sigma (\G)}$.        Therefore we have Thom collapse maps,
 
\begin{align}
\tau_\G:   \prod_{i=1}^p LM_{\ce_i}   &\xr{project}  \prod_{i=1}^p LM_{\ce_i}  / \left( (\prod_{i=1}^p LM_{\ce_i} )  -\tilde{ev}^{-1}(\eta (\Delta_\G))\right)  \cong   Map_\ce (r(\G), M)^{\tilde{ev})^*\nu (\Delta_\G)} \simeq (P_\G)^{(\tilde{ev})^*\nu (\Delta_\G)}   \quad \text{and} \notag \\
\tau_\G:  (LM)^p   &\xr{project}  (LM)^p /\left((LM)^p  - {ev}^{-1}(\eta (\Delta_\G))\right) \cong Map (r(\G), M)^{(ev)^*\nu (\Delta_\G)}
\end{align} where the targets of these maps  are the Thom spaces.  Moreover these maps are 
  compatible  in the sense that the following diagram commutes:   
$$
\begin{CD}
 \prod_{i=1}^p LM_{\ce_i}    @>\tau_\G >>  Map_\ce (r(\G), M))^{\tilde{ev}^*\nu (\Delta_\G)} @>\simeq >> (P_\G)^{(\tilde{ev})^*\nu (\Delta_\G)}   \\ 
@VVV   @VVV \\
(LM)^p @>\tau_\G >> Map (r(\G), M)^{(ev)^*\nu (\Delta_\G)}
\end{CD}
$$
 Notice that in this diagram, we again have that the  vertical maps are  homotopy equivalences. 

The bottom horizontal map, together with the Thom isomorphism, was what defined the umkehr map on the chain level in \cite{CG}
 $$
 (\rho_{in})_! : C_*(LM)^{\otimes p} \to C_{*+\chi (\G)n}(Map (r(\G), M).
 $$
 This is then compatible, up to chain homotopy,  with the Morse theoretic umkehr map,

  $$
 (\rho^{morse}_{in})_!  
  : \bigotimes_{i=1}^p C_*^{\ce_i} (LM)  \to  C_{*+\chi(\G)n}(\cmlm)
$$ defined 
to be the composition
  \begin{equation}\label{umkehr}
\begin{CD}
 (\rho^{morse}_{in})_!  :\bigotimes_{i=1}^p C_*^{\ce_i} (LM) @> h_*>\cong >  C_*( \prod_{i=1}^p LM_{\ce_i})   @>(\tau_\G)_* >> C_*(P_\G^{\nu (\Delta_\G)}) \cong   C_*( \cmlm^{\tilde \nu(\G)})  \\
&@>\cap u_\G>\cong >  C_{*+\chi (\G)n}(\cmlm)  
\end{CD}
\end{equation}
where the first map $h_*$ identifies  the Morse chain complex of the energy functional $\ce_i$ with the    cellular
complex of $LM_{\ce_i}$ which then sends it in a canonical way  to the singular chain complex $C_*(LM_{\ce_i})$.   In this diagram the symbol $``\cong "$ denotes chain homotopy equivalence, and the cap operation $\cap u_\G$ gives the Thom isomorphism on homology.     
 
We can then define the following Morse-string topology operation on the chain level,  analogous to the construction in    (\cite{CG}).
 
\begin{equation}\label{morse} 
 \mu_\G: \bigotimes_{i=1}^p C_*^{\ce_i} (LM) \xr{ (\rho^{morse}_{in})_!} C_*( \cmlm^{\tilde \nu(\G)}) \xr{ (\rho_{out})_* }  C_*(\prod_{j=1}^q LM_{\ce_{p+j}}). \end{equation}
 
 \med
 Notice that the compatibility of the umkehr map defined in \cite{CG} with the Morse umkehr map implies the following diagram of homology groups   commutes:
 
 $$
 \begin{CD}
 H_*((LM ^p))   @>(\rho^{morse}_{in})_! >> H_{*+\chi (\G)n}(\cmlm)  @>(\rho_{out})_* >> H_{*+\chi (\G)n}((LM)^q) \\
 @V=VV    @VVV  @VV=V \\
 H_*((LM ^p)) @>>(\rho_{in})_! >  H_{*+\chi (\G)n}(Map (r(\G), M))   @>>(\rho_{out})_* > H_{*+\chi (\G)n}((LM)^q) 
 \end{CD}
 $$
 In other words, we've proven the following theorem.

 \med
 
 \begin{theorem}\label{same}  For a connected surface $\Sigma$ of genus $g$, with $p$-incoming, and $q$-outgoing boundary circles, let
 $$
 \mu^{top}_{g,p+q} :  H_*((LM ^p) \to  H_{* + \chi(\Sigma)n}((LM)^q)
 $$
 be the string topology operation defined in \cite{CS} and \cite{CG}.  Then
 this operation is equal to the Morse theoretic operation
 $$
  \mu^{top}_{g,p+q} = \mu_\G
  $$
  for any connected Sullivan chord diagram $\G$ of topological type $(g; p,q)$.
  \end{theorem}
 
 \subsection{Morse theoretic string topology operations coming from families of graphs}

 In \cite{godin}, V. Godin described ``higher" string topology operations that  are  indexed  by  the homology of the moduli spaces of bordered Riemann surfaces.  In this subsection we indicate how the Morse theoretic approach to string topology described above can be adapted, using Godin's work, to yield these higher order operations.  A key ingredient in Godin's work was the generalization of the notion of a Sullivan chord diagram to a more general type of fat graph, that she called ``admissible", that had two main features: 1.  The space of admissible, marked metric fat graphs are homotopy equivalent to moduli space, and 2.  These types of graphs are sufficiently explicit so that they can be used to define the necessary umkehr maps for the definition of (higher) string topology operations. 
 These graphs were defined as follows.

\med
\begin{definition} An ``admissible" marked fat graph is one with the property that for every oriented edge $E$ that is part of an incoming boundary cycle,  its conjugate $\bar E$ (i.e the same edge with the opposite orientation) is part of an \sl outgoing \rm boundary cycle.
\end{definition}

In \cite{godin} it was proved that the space of admissible, marked fat graphs of topological type $(g, p+q)$,  $\cg_{g, p+q}$ is homotopy equivalent to the moduli space of bordered surfaces, $\cm_{g, p+q}$.   Furthermore,  if one lets $\cg_{g, p+q}(LM)$ be the space of pairs,
\begin{equation}
\cg_{g, p+q}(M) = \{(\G, \phi) \, : \, \G \in \cg_{g, p+q},\, \text{and} \,  \phi :  \G \to M  \, \text{is a continuous map} \},
\end{equation} then    $\cg^\ce_{g, p+q}(M)$ is homotopy equivalent to the space $\cm_{g, p+q}(M)$ defined in the introduction.  Furthermore the following correspondence diagram is homotopy equivalent to diagram \ref{restrict} of the introduction, and extends  diagram \ref{restrict2}:  

\begin{equation}\label{restrict3}
(LM)^q \xleftarrow{\rho_{out}} \cg_{g, p+q}(M)\xr{\rho_{in}}  \cg_{g, p+q} \times  (LM)^p .
  \end{equation}

In \cite{godin} Godin defined a generalized Pontrjagin-Thom map, which in turn defined an umkehr map $$(\rho_{in})_! : H_*( \cg_{g, p+q} \times  (LM)^p) \to H_{*+(2-2g+p+q)n}(\cg_{g, p+q}(M)).$$
 The higher string topology operations were defined as the composition
 $$
 \mu_{g, p+q} = (\rho_{out})_* \circ (\rho_{in})_!  : H_*(\cm_{g, p+q}) \otimes H_*(LM)^{\otimes p} \cong  H_*( \cg_{g, p+q} \times  (LM)^p) \to H_*((LM)^q) \cong H_*(LM)^{\otimes q}.
 $$ 
 Here we are taking homology with coefficients in an arbitrary field $k$, and the tensor products are taken over $k$. 
 
 A particular difficulty in generalizing the above Morse theoretic construction to this setting, is the technical problem that, unlike for a Sullivan chord diagram,  for a general admissible graph $\G$, the inclusion of the incoming circles,
 $$
 \alpha^- : \coprod S^1 \to \G
 $$
 (\ref{alphas}) is not an inclusion of a subcomplex (cofibration).  This will mean that our proof of the analogue of Theorem \ref{main} will not go through in this more general setting.  We get around this by considering the following larger moduli space of gradient flow surfaces defined as follows.
 
 \begin{definition}\label{cmlm1}  Let $\G$ be a marked admissible fat graph with Lagrangian labeling $\ce = \ce (\G)$.    Define the space  $\cm^{\ce}_\G(LM)_1$ to be the space of maps
$$
  \phi:  \Sigma_{\G} \to M
$$ so that the restrictions to the incoming cylinders $(-\infty, -1) \times S^1$ and the outgoing cylinders, $(+1, +\infty) \times S^1$ satisfy the gradient flow equations determined by the Lagrangian labelling.
\end{definition}

Notice that the difference between the space $\cm^{\ce}_\G(LM)$ and $\cm^{\ce}_\G(LM)_1$ is that for $\phi \in \cm^{\ce}_\G(LM)_1$ the gradient flow equations need only be  satisfied on
the cylinders $(-\infty, -1) \times S^1$  and $(+1, +\infty)\times S^1$ rather than on  entire  cylinders $(-\infty, 0) \times S^1$ and $(0, +\infty) \times S^1$ respectively.  This seemingly arbitrary distinction is important because the inclusion of each of  the incoming circles $\{-1\} \times S^1   \hk (-\infty, 0) \times S^1 \hk \Sigma_\G$ are cofibrations for all admissible fat graphs $\G$, but the inclusions of the incoming boundary circles, 
$\alpha^-_{i} : \{0\} \times S^1  \to \G \hk \Sigma_\G$ may not be cofibrations (however they would be if $\G$ were a chord diagram).  

With this technical distinction, we can now prove the following analogue of Theorem \ref{main}.

\begin{theorem}\label{main2}
  Let $\G$ be a marked admissible fat graph. 
Then the  natural map from the space of gradient flow surfaces to the continuous mapping space,
$$\cmlm_1 \to Map (\Sigma_\G, M)  $$
is a homotopy equivalence.
\end{theorem}

\begin{proof}   Since the inclusion $$\coprod_{i=1}^p \{-1\} \times S^1  \hk \coprod_{i=1}^p  (-\infty, 0)  \times S^1   \xr{\alpha^-} \Sigma_\G$$ is a cofibration,  the induced adjoint restriction map
$$
Map (\Sigma_\G, M) \to \prod_{i=1}^p LM
$$
is a fibration.  One then sees that the following 
commutative square is a pullback square of fibrations,
$$
\begin{CD}
\cmlm_1  @>>>    Map (\Sigma_\G, M)  \\
@VVV   @VVV \\
\prod_{i=1}^p LM_{\ce_i}   @>>>   (LM)^p.
\end{CD}
$$
Since, by the Palais-Smale condition, the bottom horizontal map is a homotopy equivalence, we may conclude that the top horizontal map is a homotopy equivalence.
\end{proof}

\med
We now let $\cm^\ce_{g, p+q}(LM)_1$ be the space of pairs,  
$$
\cm^\ce_{g, p+q}(LM)_1 = \{(\G, \phi): \, \G \in \cg_{g, p+q}, \, \text{and} \, \phi \in \cmlm_1\}.  
$$
Then Theorem \ref{main2}  implies that the fibration sequence  
$$
\cmlm_1  \to \cm^\ce_{g, p+q}(LM)_1 \to  \cg_{g, p+q}
$$
is homotopic to the fibration sequence
$$
Map (\Sigma_G, M) \to \cg_{g, p+q}(M) \to \cg_{g, p+q},
$$
which, by Godin's result \cite{godin}   is in turn homotopic to the fibration sequence,
$$
Map (\Sigma_G, M) \to \cm_{g, p+q}(M) \to \cm_{g, p+q}.
$$
We therefore have a commutative diagram
$$
\begin{CD}
\prod_{i=1}^q LM_{\ce_j} @<\rho_{out} <<  \cm^\ce_{g, p+q}(LM)_1  @>\rho_{in} >>   \cg_{g, p+q} \times \prod_{i=1}^p LM_{\ce_i} \\
@V\simeq VV     @VV\simeq V  @VVV  \\
(LM)^q    @<<\rho_{out} <   \cg_{g, p+q}(M) @>>\rho_{in}>  \cg_{g, p+q} \times (LM)^p 
\end{CD}
$$
where the horizontal maps are homotopy equivalences. Using these equivalences and Godin's    construction one obtains higher order operations,   
 $$
 \mu^{morse} : H_*(\cm_{g, p+q}) \otimes \bigotimes_{i=1}^p  H_*(LM_{\ce_i})     \to H_{* }(\cm^\ce_{g, p+q}(LM)_1)   \to \bigotimes_{i=1}^{p} H_*(LM_{\ce_i}).
 $$


\section{String operations by counting gradient flow lines}
 
 In this section we give a more analytical description of the Morse theoretic string topology operations, via the counting of zero dimensional moduli spaces of gradient trajectories.  Such an operation was constructed in \cite{abboschwarz2} corresponding to the figure 8 graph.  Our comparison of these operations with the ones constructed in section 1, will imply that this figure 8 product on the Morse homology of $LM$ indeed corresponds to the Chas-Sullivan loop product.  Combining this with the theorem of \cite{abboschwarz2} giving a ring isomorphism between the Floer homology of the cotangent bundle $HF_*(T^*M)$ with the Morse homology $H_*(LM)$  implies that the pair of pants product in Floer homology corresponds to the Chas-Sullivan product in $H_*(LM)$.
 
We continue to consider a Lagrangian, energy functional, and metric satisfying the conditions described in section 1.

Let $\G$ be a Sullivan chord diagram of type $(g; p,q)$ and $r(\G)$ the corresponding reduced chord diagram with parametrizations $\alpha^-\colon \coprod_q S^1 \la r(\G)$ for the incoming cycles and $\alpha^+\colon \coprod_p S^1\la r(\G)$ for the outgoing cycles, where we reparametrize $\alpha^\pm$ such that $S^1=\br/\bz$ is the standard circle.

The space $Map\big(r(\G),M\big)$ is endowed with the structure of a Hilbert manifold, the topology given by edgewise Sobolev $W^{1,2}$-maps. The parametrizations $\alpha^\pm$ induce embeddings
$$
   (LM)^q \stackrel{r_{out}}{\hookleftarrow} Map\big(r(\G),M\big) \stackrel{r_{in}}{\hookrightarrow} (LM)^p
$$
by $r_{in}(c)=c\circ\alpha^-$ and $r_{out}(c)=c\circ\alpha^+$.

Recall the following pullback square of fibre bundles:

\begin{equation} 
  \begin{CD}
  Map(r(\G), M)   @>r_{in} >>    (LM)^p \\
  @Vev_\G VV    @VVev_\G V  \\
  M^{\sigma (\G)}  @>>\Delta_\G >  M^{v(\G)}
  \end{CD}
\end{equation}
Since this is a pullback of smooth bundles, $ Map(r(\G), M)$ has the structure of a codimension $(-\chi (\G)\cdot n)$-submanifold of $(LM)^p$ with coorientation induced by the embedding $\Delta_\G$, as $M$ is assumed oriented.   Recall that in the definition of the loop product, one uses the figure 8 graph for  $\G $,  which is already a reduced chord diagram.  Similarly, the ``little cacti" diagrams (see \cite{CJ}, \cite{cohenvoronov}, \cite{voronov}) are also reduced Sullivan chord diagrams.  These are the diagrams used to describe the ``BV"-structure in string topology.   
  
For our purposes it is not necessary to describe the analogous structure of a cooriented embedding for $r_{out}$. However, we require that the following transversality conditions hold. 
 
\med
\noindent \bf Transversality Condition 1. \sl For any collection of critical points  \sl   $a_i\in\cp(L_i)$, $i=1,\ldots,p$ the embedding  $r_{in}\colon Map\big(r(\G),M\big) \hk (LM)^p$ is transverse to
$W^u(a_1) \times \cdots \times W^u(a_p) \hk LM \times \cdots \times LM$.

\med
\rm and similarly
 \med
 
 \noindent \bf Transversality Condition 2. \sl For any collection of critical points  \sl   $a_{p+i}\in\cp(L_{p+i})$, $i=1,\ldots,q$, the embedding  $r_{out}\colon Map\big(r(\G),M\big) \hk (LM)^q$ is transverse to
$W^s(a_{p+1}) \times \cdots \times W^s(a_{p+q}) \hk LM \times \cdots \times LM$.

\rm
 
The first transversality condition implies that the intersection
$$
  \cm_{r(\G)}(LM;a_1,\cdots,a_p)\,=\,r_{in}\big(Map(r(\G),M)\big)\pitchfork W^u(a_1)\times\ldots\times W^u(a_p)
  $$
is an oriented, finite-dimensional submanifold of $(LM)^p$ of dimension $\sum_{i=1}^p Ind(a_i)+\chi(\G)\cdot n$, and analogously for Condition 2.

However, in order to construct the Morse-theoretical description of $\mu_\G$, we need a stronger condition than Transversality Condition 2. Condition 2 will be used within the proof of the following theorem.

 \med
 
 \noindent \bf Transversality Condition 3. \sl For any collection of $p+q$ critical points  \sl   $a_i\in\cp(L_i)$, $i=1,\ldots,p+q$, the restriction   
$$
    r_{out}{}_{|\cm_{r(\G)}}\colon  \cm_{r(\G)}(LM;a_1,\cdots, a_p) \hk (LM)^q 
$$    
  is transverse to $W^s(a_{p+1}) \times \cdots \times W^s(a_{p+q}) \hk LM \times \cdots \times LM$.

\rm
  
 \med
 We will discuss these tranversality conditions below.
Assuming them for now,    we have for 
 $\vec{a} = (a_1, \cdots , a_p, a_{p+1}, \cdots, a_{p+q})\in\prod_{i=1}^{p+q}\cp(L_i)$   the following immediate result:

 \begin{proposition}\label{dimension} 
 The space
 $$
    \cm_{r(\G)}(LM;\vec{a}):=r_{out}\big(\cm_{r(\G)}(LM;a_1,\cdots,a_p)\big) \pitchfork
    W^s(a_{p+1})\times\cdots\times W^s(a_{p+q})\,\subset\, (LM)^q
    $$
is a   smooth orientable  manifold of dimension 
$$
   \dim\cm_{r(\G)}(LM;\vec{a})\,=\,
   \sum_{i=1}^p Ind(a_i) - \sum_{j=p+1}^{p+q} Ind(a_j) +\chi (\G)\cdot n
   $$ 
homeomorphic to
$$
   \big\{\,\phi\in\cm^\ce_{r(\G)}(LM)\,\big|\,
   \lim_{t\to-\infty}\phi_i(t,\cdot)=a_i,\,i=1,\cdots, p,\;
   \lim_{t\to+\infty}\phi_j(t,\cdot)=a_j,\,j=p+1,\cdots, p+q,\big\}\,.
   $$
An orientation of this manifold is induced by orientations of the tangent spaces of the unstable manifolds $T_{a_i}W^u(a_i)$ of the critical points $a_1,\cdots,a_{p+q}$.
 \end{proposition}
 
 \med
 In the case when $\vec{a}$ is a collection of critical points so that the dimension of the manifold $\cm_{r(\G)}(LM;\vec{a})$ is zero, standard considerations imply that it is compact  (see for example \cite{abboschwarz2}).
 We may then  define an  operation on the Morse chain complex (\ref{complex}).
 
 \begin{definition}  Define the operation $$\nu_\G : C^\ce_*((LM^p)   \to C^\ce_{*+\chi(\G)n}((LM)^q)
 $$
 by 
 \begin{equation*}
 \nu_\G([a_1] \otimes \cdots \otimes [a_p]) =  
 \hspace{-10em} \sum_{\hspace{10em}dim \, \cm_{r(\G)}(LM; a_1, \cdots a_p, a_{p+1}, \cdots a_{p+q})=0}  \hspace{-11em}\#\cm_{r(\G)}(LM; a_1, \cdots a_p, a_{p+1}, \cdots a_{p+q}) \, [a_{p+1}] \otimes \cdots \otimes  [a_{p+q}]  
 \end{equation*}
 Here $\#\cm_{r(\G)}(LM;\vec{a})$ is the  oriented count of the number of points in this  zero dimensional compact manifold.   
  \end{definition}
 
Concerning the previous transversality conditions, they are fulfilled for generic choices of Riemannian metrics on the Hilbert manifold $LM$, provided that the Lagrangians $L_i$ have been chosen suitably. Namely, the solution sets of gradient flow trajectories corresponding to the intersections in question should not contain constant solutions. For example, in Condition 1, the $p$-tuple $(a_1,\dots,a_p)\in \Crit\ce_{L_1}\times\ldots\times\Crit\ce_{L_p}$ should not be contained in $Map\big(r(\G),M\big)$. Moreover, the proof of the generic existence of such metrics uses the theorem of Sard-Smale. All of the above intersection problems are Fredholm problems, but in general we need to consider not only Fredholm indices up to 1, but also higher, as e.g. in Transversality Condition 1. This is the place where we have to assume accordingly high differentiability of our energy functions $\ce_{L_i}$, $i=1,\ldots,p+q$.

Assume for example that $\cb$ is a suitably defined separable Banach manifold consisting of admissible Riemannian metrics on $LM$. Then we consider the infinite-dimensional Banach manifold
$$
   \cw^u(a_1,\dots,a_p) = \big\{\,(g,c_1,\ldots,c_p)\in\cb\times LM^p\,|\,c_i\in W^u_g(a_i),\,i=1,\ldots,p\,\}
   $$
where $W^u_g(a_i)$ is the unstable manifold for the negative gradient flow of $\ce_{L_i}$ determined by the Riemannian metric $g$. It not hard to show that for a sufficiently rich set $\cb$ of variations of the Riemannian metric on $LM$, each projection $\cw^u(a_1,\ldots,a_p)-\{(a_1,\ldots,a_p)\}\to LM$, $(g,\vec{c})\mapsto c_i$ is a submersion onto its image away from the critical point. This, together with the assumed smoothness is the main ingredient in the application of the Sard-Smale theorem in order to prove that the Transversality Condition 1 is generically fulfilled. Similarly, we also obtain the other transversality conditions as generically satisfied. For more details on this transversality analysis we refer to \cite{abbomajer06,abboschwarz}.
 
 \med
 The following is now the main theorem of this section.
 
 \med
 \begin{theorem}\label{same}
Under the above assumptions on the Lagrangians and the metric, and assuming transversality conditions 1 and 3,  the operation $\nu_\G$ is a chain map, and in homology it gives the string topology operation
$$
\nu_\G = \mu_\G : H_*((LM)^{p})   \to H_{*+\chi(\G)n}((LM)^q).
$$
\end{theorem}
 
For the proof, we factorize $\nu_\G$ in close analogy to $\mu_\G$ above, into
\begin{equation}\begin{split}
  &\nu_\G\,=\, (r_{out})_\ast \circ (r_{in})_!,\quad\text{where}\\
  &(r_{in})_!\colon  H^\ce_\ast\big( (LM)^p\big) \la H^\cf_{\ast+\chi(\G)\cdot n}\big (Map(r(\G),M)\big),\quad\text{and}\\
  &(r_{out})_\ast\colon H^\cf_\ast\big( Map(r(\G),M)\big) \la H^\ce_\ast\big( (LM)^q\big)
  \end{split}
\end{equation}
are naturally isomorphic to $(\rho_{in})_!$ and $(\rho_{out})_\ast$ when we compare Morse homology with standard homology.

Here we consider an auxilary smooth Morse function $\cf$ on the Hilbert manifold $Map(r(\G),M)$. That is, $\cf$ satisfies the Palais-Smale property
with a complete negative gradient flow for a complete Riemannian metric, it is bounded below, and all critical points $b\in\Crit \cf$ are non-degenerate and of finite Morse index. Choosing a generic Riemannian metric on $Map(r(\G),M)$ satisfying Morse-Smale transversality with respect to $\cf$ for relative Morse index up to 2, the Morse complex
$\big( C^\cf_\ast(Map(r(\G),M)),\partial\big)$ is well-defined and its Morse homology naturally isomorphic to standard homology $H_\ast\big(Map(r(\G),M)\big)$, see e.g. \cite{abbomajer,abboschwarz2,Schwarz-mono}.

Essentially, up to a small perturbation, we can take for $\cf$ the restricted energy functional $r_{in}^\ast\ce^{\otimes p}=\ce^{\otimes p}{|r_{in}(Map(r(\G),M))}$, for $\ce^{\otimes p}\colon (LM)^p\to\br$ given by $\ce^{\otimes p}(c_1,\ldots,c_p)=\ce_{L_1}(c_1)+\ldots +\ce_{L_p}(c_p)$.

\med

We now focus on the embedding $r_{out}\colon Map(r(\G),M) \hookrightarrow (LM)^q$. In addition to Transversality Condition 2, we assume

 \med
 
 \noindent \bf Generic Condition 4. \sl 
  $r_{out}$ maps no critical point of $\cf$ to a critical point of $\ce$.

\rm
\med
Given $b\in\Crit\cf$ and $a_{p+i}\in\cp(L_{p+i})$, $i=1,\ldots, q$, we define
$$
   \cm^{out}_{r(\G)}(b,a_{p+1},\cdots,a_{p+q})\,=\,
   r_{out}\big( W^u_\cf(b)\big) \cap \big( W^s_\ce(a_{p+1})\times\cdots\times W^s_\ce(a_{p+q})\big),
$$
where $W^u_\cf(b)$ is the unstable manifold for the negative gradient flow of $\cf$. Condition 4  allows us to find a generic metric on $Map(r(\G),M)$ such that $W^u_\cf(b)$ intersects the submanifold $r_{out}^{-1}\big( W^s(a_{p+1})\times\cdots\times W^s(a_{p+q})\big)$ transversely.  Altogether  we obtain $\cm^{out}_{r(\G)}(b,a_{p+1},\cdots,a_{p+q})$ as a manifold of dimension $ind(b) -\sum_{j=1}^q ind(a_{p+j})$ with orientation induced by the orientations  of the unstable manifolds  $W^u(b)$, $W^u(a_{p+1}),\cdots, W^u(a_{p+q})$.

We define
\begin{equation}\begin{split}
  &(r_{out})_\ast\colon C^\cf_\ast \big( Map(r(\G),M)\big) \la C^\ce_\ast\big( (LM)^q\big),\\
  &(r_{out})_\ast ([b]) \,=\, \hspace{-6em}  \sum_{\hspace{6em}  \dim\cm^{out}_{r(\G)}(b,a_{p+1},\cdots,a_{p+q})=0}
     \hspace{-7em} \#\cm^{out}_{r(\G)}(b,a_{p+1},\cdots,a_{p+q})\,[a_{p+1}]\otimes\cdots\otimes [a_{p+q}]
     \end{split}
\end{equation}
and have the following result. 

\begin{proposition}
  $(r_{out})_\ast\colon C_\ast^\cf\big( Map(r(\G),M)\big)\la C_\ast^\ce\big( (LM)^q\big)$ is a chain map and the natural isomorphism to standard homology intertwines $(r_{out})_\ast$ with $(\rho_{out})_\ast$.  That is, the following diagram commutes:  
 $$ \begin{CD}
     H_\ast^\cf\big( Map(r(\G),M)\big)  @>(r_{out})_\ast>> H_\ast^\ce\big( (LM)^q\big) \\
	@VV\cong V  @VV\cong V\\
	H_\ast \big( Map(r(\G),M)\big)  @>(\rho_{out})_\ast>> H_\ast \big( (LM)^q\big) \\
  \end{CD}
  $$
  commutes.
\end{proposition}

This construction describes the functoriality for Morse homology. More details can be found in \cite{abboschwarz2} and \cite{Schwarz-mono}.

\med

We now construct the Morse-theoretical version $(r_{in})_!$ of the umkehr map. Recall the map 
$$r_{in}\colon Map(r(\G),M)\hookrightarrow (LM)^p$$ which is a proper embedding of Hilbert manifolds, with finite codimension and cooriented. In addition to Transversality Condition 1, we assume in correspondence to Condition 4
 \med
 
 \noindent \bf Generic Condition 5. \sl 
$r_{in}$ maps no critical point of $\cf$ to a critical point of $\ce$.

\rm
\med
Hence, we find a generic metric on $Map(r(\G),M)$ such that for any $b\in\Crit\cf$ its stable manifold $W^s_\cf(b)$ is transverse to $\big( W^u_\ce(a_1)\times\cdots\times W^u_\ce(a_p)\big)\cap r_{in}\big( Map(r(\G),M\big)$ for all $a_i\in\cp(L_i)$, $i=1,\ldots,p$. We obtain the manifold
$$
   \cm^{in}_{r(\G)}(a_1,\cdots,a_p,b)\,=\,\big( W^u_\ce(a_1)\times\cdots\times W^u_\ce(a_p)\big)\cap W^s_\cf(b)
   $$
of dimension $\sum_{i=1}^p ind(a_i) - ind(b) -\chi(\G)\cdot n$ and induced orientation from the orientations of $W^u_\ce(a_1),\cdots, W^u_\ce(a_p)$, coorientation of $W^s_\cf(b)$ within $Map(r(\G),M)$ by the orientation of $W^u_\cf(b)$, and the coorientation of $r_{in}\big( Map(r(\G),M)\big)$ within $(LM)^p$.

We define
\begin{equation}  \begin{split}
  &(r_{in})_!\colon C^\ce_\ast\big( (LM)^p\big) \la C^\cf_{\ast+\chi(\G)\cdot n} \big( Map(r(\G),M)\big),\\
  &(r_{in})_!\big( [a_1]\otimes\cdots\otimes [a_p]\big)\,=\hspace{-6em}  \sum_{\hspace{6em}  \dim\cm^{in}_{r(\G)}(a_{1},\cdots,a_{p},b)=0}
     \hspace{-6em} \#\cm^{in}_{r(\G)}(a_{1},\cdots,a_{p},b)\,[b]\,.
   \end{split}
\end{equation}   
By the usual arguments of Morse homology, this is a chain map. Note that the compactness of $\cm^{in}_{r(\G)}$ also requires the properness of the embedding $r_{in}$, such that $r_{in}^{-1}\big(W^u(a_1)\times\cdots\times W^u(a_p)\big)$ is relatively compact in $Map(r(\G),M)$. 

We have
\begin{proposition}
    The chain map $(r_{in})_!\colon C^\ce_\ast\big( (LM)^p\big) \la C^\cf_{\ast+\chi(\G)\cdot n}\big( Map(r(\G),M)\big)$ induces an umkehr map on Morse homology compatible with the umkehr map $(\rho_{in})_!$ under the natural isomorphism, i.e.
     $$ \begin{CD}
     H_\ast^\ce\big( (LM)^p\big )  @>(r_{in})_!>> H_{\ast+\chi(\G)\cdot n}^\cf\big(Map(r(\G),M)\big) \\
	@VV\cong V  @VV\cong V\\
     H_\ast\big( (LM)^p\big )  @>(\rho_{in})_!>> H_{\ast+\chi(\G)\cdot n}\big(Map(r(\G),M)\big)
  \end{CD}
  $$
  commutes.
\end{proposition}

Before we give  a proof of this by a Morse-theoretical description of the Thom isomorphism, we conclude the proof of Theorem \ref{same}. Via a standard gluing argument, $(r_{out})_\ast\circ (r_{in})_!$ on Morse chain level is equal to counting
\begin{equation*}\begin{split}
   &\Theta\in r_{in}^{-1}\big( W^u(a_1)\times\cdots\times W^u(a_p)\big)\,\subset\, Map(r(\G),M),\quad \text{such that}\\
   &\phi^R_\cf(\Theta)\in r_{out}^{-1}\big( W^s(a_{p+1})\times\cdots\times W^s(a_{p+q})\big),
   \end{split}
\end{equation*}
for $R>0$ fixed and sufficiently large, where $t\mapsto\phi^t_\cf(\Theta)$ is a flow line for the negative gradient flow of $\cf$ on $Map(r(\G),M)$.   

By homotoping $R$ to $0$ we establish a cobordism to the previous solution space $\cm_{r(\G)}(LM;\vec{a})$ which gives rise to a chain homotopy operator, proving
$$
  \nu_\G\,\simeq\, (r_{out})_\ast\circ (r_{in})_!
  $$
on chain level. Hence, via the natural functor to standard homology, we have on homology level
$$
   \nu_\G\,=\,   (r_{out})_\ast\circ (r_{in})_! \,=\, (\rho_{out})_\ast\circ (\rho_{in})_! \,=\, \mu_\G,
   $$
proving Theorem \ref{same}.

 \section{Thom-Isomorphism, Euler-Class and the Umkehr Map via Morse Homology}

We will now give a Morse-theoretical construction of the umkehr map which is based on a Morse-theoretical construction of the Thom isomorphism. We need to give this construction in the infinite-dimensional setting, at least sufficient for the case of the loop space $LM$ with its $W^{1,2}$-Hilbert manifold structure and the energy functional $\ce$. In fact, Morse homology can be defined for a much larger class of infinite-dimensional settings. For more details on the difference between the finite and the infinite dimensional case and for the more general setting of the latter we refer to \cite{abbomajer,abbomajer05}.

\subsection{Preliminaries} Let $X$ be a smooth paracompact Hilbert manifold with a complete Riemannian metric, and let $f\colon
X\to\R$ be a $C^2$ Morse function satisfying the Palais-Smale condition. Moreover, for our purpose ($X=LM$, $f=\ce$) we assume that $f$ is bounded below and all critical points are of finite Morse index.

Pick a generic Riemannian metric $g$, s.t. $(f,g)$ is a Morse-Smale-pair \footnote{Transversality is sufficient up to index difference 2. See 0.5 in \cite{abbomajer05} for the precise details for genericity here.} with a complete negative gradient flow. For
$x,y\in\Crit f$ with $i(x)-i(y)=1$ set
$$
   \langle x,y\rangle\,=\,\#_{\mathsf{alg}}\big\{\,\gamma\colon \R\to
   X\,|\, \dot\gamma+\nabla f(\gamma)=0,\;
   \gamma(-\infty)=x,\,\gamma(\infty)=y\,\big\}\big/\R\,.
$$
Here, $\#_{\mathsf{alg}}$ refers to counting with orientations
obtained from the concept of coherent orientations, that is arbitary
orientations of all unstable manifolds and induced coorientations of
all stable ones.
We obtain the Morse homology of $f$ from
\begin{equation*}\begin{split}
  &C_\ast(f)\,=\Z\otimes\Crit_\ast f,\qquad\text{(finitely generated)}\\
  &\partial\colon C_\ast(f)\to C_{\ast-1}(f),\quad \partial x=
  \sum_{i(y)=i(x)-1}\langle x,y\rangle y, 
\end{split}
\end{equation*}
and the Morse cohomology from
\begin{equation*}\begin{split}
    &C^\ast(f)=\Z^{\Crit_\ast f},\qquad\text{(not necessarily finitely generated !)}\\
  &\delta\colon C^\ast(f)\to C^{\ast+1}(f),\quad
  (\delta\phi)(x)=\sum_{i(y)=i(x)-1} \langle x,y\rangle \phi(y),\\
  &\text{i.e. }\quad  \big(C^\ast(f),\delta\big)\,=\,\operatorname{Hom}\big((C_\ast(f),\partial),\Z\big)\,.
\end{split}
\end{equation*}
By identifying $\delta_x\in\Z^{\Crit f}$ with $x\in\Crit f$ we see
that $\delta$ is alternatively defined by counting positive gradient
flow lines for $f$.

Note that from the lower boundedness of $f$, the completeness of the negative gradient flow and the Palais-Smale property, we obtain that $H_0(f)\cong\bz\cong H^0(f)$ if $X$ is connected. Hence we have a generator $1\in H_0(f)$ represented by a single critical point of index 0 and $1\in H^0(f)$ represented by $\phi\in\bz^{\Crit_0 f}$, $\phi\equiv 1$. 

\subsection{Relative Cohomology}
We now recall the Morse-theoretical definition of relative homology
and cohomology, see e.g. \cite{Schwarz-mono}. Let $A\subset X$ be an open submanifold with smooth
boundary $\partial A$, and we assume that the above Morse function $f$ on $X$  is in addition such that $\nabla
f\pitchfork\partial A$, and the gradient $\nabla f$ is pointing out of
$A$. This implies that
$$
   \big(C_\ast(f_{|A}),\partial\big)\quad
   \text{is a subcomplex of }\quad \big(C_\ast(f),\partial\big),
$$
and we have the exact sequence of chain complexes
$$
   0\to C_\ast(f_{|A})\stackrel{i}{\longrightarrow} C_\ast(f)
   \stackrel{j}{\longrightarrow} C_\ast\big(f;X,A):=C_\ast(f)/C_\ast(f_{|A})\to 0,
$$
inducing the long exact sequence of homology.

For Morse cohomology we have dually
$$
 0\to C^\ast(f;X,A)\stackrel{j^\ast}{\longrightarrow}
 C^\ast(f;X)\stackrel{i^\ast}{\longrightarrow} C^\ast(f_{|A};A)\to
 0\,.
$$
Namely, we have
$$
   i^\ast(\sum_{x\in\Crit f} a_x x)\,=\,\sum_{x\in\Crit f\cap A} a_x x,
$$
since $x\in \Crit f\cap (X\setminus A)$ implies 
$\delta x\in \Z^{\Crit  f\cap (X\setminus A)}$, hence $i^\ast$ is a
cochain complex morphism. Also, we set
$$
  C^\ast(f; X,A)\,=\,\Z^{\Crit f\cap(X\setminus A)}
$$
which by the same argument of $\nabla f$ pointing outwards along
$\partial A$ turns $j^\ast$ into a sub-cochain-complex inclusion. Note
that here the obvious excision principle is used in the notation
$$
   C^\ast(f;X,A)\,=\,C^\ast(f;X\setminus A,\partial A)\,.
$$

\subsubsection*{Example}
Let us consider the following very simple example, which
describes the main idea used for the following Thom
isomorphism and which illustrates the importance of changing Morse functions in Morse homology. Let $q\colon \R^n\to\R$ be the standard positive
quadratic form with its unique critical point of index $0$ in the
origin. Choose another coercive Morse function $\tilde q\colon \R^n\to
\R$  such that
$$
\tilde q(x)=\begin{cases}
  -q(x),\,&|x|\leq 1,\\
   q(x),\,&|x|\geq 2\,.
 \end{cases}
 $$
Obviously we have for the unit disk $D_1$ and its boundary sphere
$S_1$
$$
   H^\ast(-q;D_1,S_1)\,=\,H^\ast(\tilde q;D_1,S_1)\,.
$$
The above long exact cohomology sequence and obvious identifications
and homotopy invariance give
$$
\begin{CD}
   0@>{}>> H^{n-1}(\tilde q;\R^n\setminus D_1)@>\delta^\ast>>
   H^n(\tilde q;D_1,S_1)@>j^\ast>> H^n(\tilde q;\R^n)\\
   &&&&@| @V\cong VV\\
   &&&& H^n(-q;D_1,S_1) && H^n(q;\R^n)\\
   &&&&@V\cong VV@|\\
   &&&&\Z&&\{0\}\\
 \end{CD}
$$
This leads to $H^{n-1}(\tilde q;\R^n\setminus D_1)\cong\Z$. In fact
$\tilde q$ encodes the Morse cohomology of the $(n-1)$- sphere.
 
\subsubsection*{Functoriality for proper embeddings}
Consider now a proper embedding of finite and positive codimension of a
submanifold $e\colon P\hookrightarrow X$, and let $k\colon
P\to\R$ be a Morse function on $P$ of the same type as $f$, i.e. satisfying the Palais-Smale property, bounded below and only with critical points of finite index. We define $e_\ast
\colon C_\ast(k)\to C_\ast(f)$ as follows. Consider $p\in \Crit k\subset P$
and $x\in\Crit f\subset X$. We require generic metrics on $P$ and $X$
such that $W^u(p;k)$ as a submanifold of $X$ and $W^s(x;f)$ intersect
transversely in $X$. This can always generically be achieved for sufficiently bounded relative index in the case of critical points of $f$ in the complement of $P$. In the case, where we have to allow critical points of $f$ on $P$ we need to verify that the transverse intersection is already satisfied automatically.

 In such a transverse situation we define
$$
  n(p,x)\,=\,\#_{\mathsf{alg}}\big(W^u(p,k)\pitchfork_X W^s(x;f)\big)
$$
if $i(p)=i(x)$. Note that, for this relative Morse index, the
intersection is $0$-dimensional. We then define
$$
   e_\ast(p)=\sum_{i(x)=i(p)} n(p,x) x\,,
$$
and we easily see that
$\partial_f\circ e_\ast=e_\ast\circ\partial_k$.
Analogously, we define the pull-back homomorphism
$$
   e^\ast\colon C^\ast(f;X)\to C^\ast(k;P)\,.
$$

Consider now the special case where $\dim P=\ind(x)=l$ and 
$m\in \Crit f\cap P$, $\ind (p)=l$ for $p\in\Crit k$. If $T_x
W^s(x,f)=\operatorname{Eig}^+D^2 f(x)$ is transverse to $P$, then the
above transversality for $n(p,x)$ is automatically satisfied.

\subsubsection*{Construction of the Thom isomorphism}
Let us now consider a smooth vector bundle $\pi\colon E\to X$ of finite rank
$r$ endowed with an arbitrary Riemannian metric. Let $q\colon E\to
[0,\infty)$ be the associated positive quadratic form, and consider
the disk and sphere bundle
\begin{equation*}\begin{split}
  &D(E)\,=\,\big\{\,(x,v)\in E_x\,|\,q(v)\leq 1\,\big\},\\
  &S(E)\,=\,\big\{\,(x,v)\in E_x\,|\,q(v) = 1\,\big\}\,=\,\partial
  D(E)\,.\\
\end{split}
\end{equation*}
Thus $f_{-q}:=\pi^\ast f-q$ is an admissible relative Morse function
for the pair $(E,E\setminus D(E))$. If we extend $f_{-q}|D(E)$ to $\tilde
f_q$ outside of $D(E)$ such that $\tilde f_q$ is a Morse function and
$$
   \tilde f_q(x,v)=f(x)+q(v)\quad\text{for }\;q(v)\geq 2,
$$
we have the canonical identification as in the above example
$$
   H^\ast(f_{-q};D(E),S(E))\,=\,H^\ast(\tilde f_q;D(E),S(E))\,=\,
   H^\ast(\tilde f_q;E,E\setminus D(E))
$$
and the exact sequence
\begin{equation}\label{seq1}
  \ldots\to H^{\ast-1}(\tilde f_q;E\setminus D(E))\,\to\,
  H^\ast(\tilde f_q;D(E),S(E))\stackrel{j^\ast}{\longrightarrow}
  H^\ast(\tilde f_q;E)\,\cong\,H^\ast(f_q;E),
\end{equation}
with $f_q=\pi^\ast f+q$. Moreover, there is the canonical isomorphism
\begin{equation}\label{iso1}
   \pi^\ast\colon H^\ast(f;X)\stackrel{\cong}{\longrightarrow}
   H^\ast(f_q;E)
 \end{equation}
 induced from
$$
   \Crit_\ast f\,=\,\Crit_\ast f_q
$$
by identifying $X$ with the zero section of $E$.
\begin{proposition}\label{prop Morse-Thom} If $E$ is an orientable bundle, then $\Crit_\ast
  f\,=\,\Crit_{\ast+r}f_{-q}$ induces an
  isomorphism 
$$
T^\ast\colon H^\ast(f;X)\stackrel{\cong}{\longrightarrow}
H^{\ast+r}(f_{-q};D(E),S(E))
$$
and the element $u=T^\ast(1)\in H^r(f_{-q})$ satisfies
$\varphi_x^\ast(u)=u_o$ where $\varphi_x\colon
(D^r,S^{r-1})\hookrightarrow (D(E),S(E))$ is the fibre inclusion over
any $x\in X$ and $u_o\in H^r(D^r,S^{r-1})$ is the generator compatible
with the orientation of $E$. Moreover, the same identification of critical points induces the dual isomorphism $T_\ast\colon H_\ast(f_{-q};D(E),S(E))\stackrel{\cong}{\longrightarrow} H_{\ast-r}(f;X)$.
\end{proposition}
\begin{corollary} $u\in H^r(f_{-q};D(E),S(E))=H^r(\tilde f_q;D(E),S(E))$ is
  the Thom class of $E$, and $(\pi^{\ast})^{-1}\circ
  j^\ast(u)=:e(E)\in H^r(X)$ is the Euler class of $E$.
\end{corollary}
This follows from the proposition using (\ref{seq1}) and (\ref{iso1}).

We now prove the proposition. 
\begin{proof} 
Without loss of generality, $f$ on $X$ has a unique minimum in
$x_o$ and we can choose $x=x_o$ for the fibre inclusion. Hence $T^\ast(1)$ is represented by $\{(x_0,0)\}\in
H^r(f_{-q};D(E),S(E))$. Obviously,
$\operatorname{Eig}^+\big((x_o,0);f_{-q}\big)\cong\operatorname{Eig}^+(x_o;f)$
is transverse to the fibre $E_{x_o}$. Moreover, we have\\
$W^s\big((x_o,0);f_{-q}{}{|E_{x_o}}\big)=\{0\}$ and
$W^u(0;-q)=T_0E_{x_o}$. Hence, we have transverse intersection within
$E$ and we see $\varphi_{x_o}^\ast(x_o,0)=0$. Obviously, $\{0\}=u_o\in
H^r(-q;D^r,S^{r-1})$ is the generator.

Using the canonical identification $\Crit_\ast f=\Crit_{\ast+r}f_{-q}$ via $x=(x,0)$, $W^s(x;f)=W^s((x,0);f_{-q})$ for all $x\in\Crit_\ast f$, and using the fact that each oriented unstable manifold $W^u(x;f)$ together with the orientation of $E$ gives an orientation for $W^u((x,0);f_{-q})$, implies that we have
$$
   \#_{\mathsf{alg}}\big( W^u(x;f)\pitchfork W^s(x';f)\big)\,=\,
   \#_{\mathsf{alg}}\big( W^u((x,0);f_{-q})\pitchfork
   W^s((x',0);f_{-q})\big)\,.  
$$
Since there are no critical points for $f_{_q}$ off the zero section
we have canonical chain and cochain complex isomorphisms
\begin{equation*}\begin{split}
   T_\bullet\colon &C_{\ast+r}\big(f_{-q};D(E),S(E)\big)
   \stackrel{\cong}{\longrightarrow} C_\ast(f;X),\\
   T^\bullet\colon &C^\ast(f;X) \stackrel{\cong}{\longrightarrow}
   C^{\ast+r}\big(f_{-q};D(E),S(E)\big) ,\\
 \end{split}
\end{equation*} 
inducing the Thom isomorphisms $T_\ast$ and $T^\ast$ on homology
respectively cohomology level.
\end{proof}

\subsection{Umkehr map for proper embeddings}
We will now give another, Morse-theoretical description for the umkehr
map.

Let $e\colon P\hookrightarrow X$ be again a proper embedding of
finite positive codimension with the additional assumption of coorientation,
i.e. the normal bundle $\nu_e$ is oriented. Consider again
Morse functions as above, $k\colon P\to \R$ and $f\colon X\to \R$ with $\Crit
f\cap P=\emptyset$. Then, for generic metrics on $P$ and $X$ we have
transverse intersections of the unstable manifold of $m\in \Crit f$ in
$X$ and the stable manifold of $p\in \Crit k$ in $P$,
$$
   W^u(m;f,X)\pitchfork W^s(p;k,P)\,.
$$
Note that, this in particular requires the tranverse intersection of
$W^u(m;f,X)$ with $P$. The coherent orientation condition of Morse
homology is guaranteed by the assumption of coorientation of $e(P)$ in
$X$. Hence, we obtain a well-defined integer
$$
   n(m,p)\,=\,\#_{\mathsf{alg}}\big(W^u(m;f)\pitchfork W^s(p;k)\big)
$$
if $i(m)=i(p)+r$, where $r$ is the codimension of $P$.
\begin{proposition} \label{prop umkehr}
The associated Morse chain morphism
$$
  e_\bullet\colon C_\ast(f;X)\,\to\, C_{\ast-r}(k;P),\quad
  m\mapsto \sum_{i(p)=i(m)-r} n(m,p)p,
$$
induces the umkehr map
$$
  e_{!}\colon H_\ast(f;X)\,\to\,H_{\ast-r}(k;P)
$$
on the level of Morse homology.
\end{proposition}
\begin{proof}
Again, it is a standard Morse homology argument to see that the above
chain level map $e_\bullet$ commutes with the respective boundary
operators, $\partial_k\circ e_\bullet=e_\bullet\circ\partial_f$.

Consider the tubular neighbourhood $\eta_e$ of $P$ and let $\tilde
k_q\in C^\infty(X,\R)$ be a Morse function on $X$ whose
restriction to $\eta_e$ equals $k_{-q}=\pi^\ast k-q$ as above, after identifying
$\eta_e$ with an open subset of the normal bundle $\nu_e$. Using the
orientation of $\nu_e$, we can identify
$$
   C_\ast(\tilde k_q;X,X\setminus\eta_e)\,\cong\,C_{\ast-r}(k;P)
$$
precisely as in the proof of Proposition \ref{prop Morse-Thom}.
Moreover, using $j$ from the short exact sequence on chain level, we have
\begin{equation}\label{eq collapse}
  C_\ast(\tilde k_q;X)\stackrel{j}{\longrightarrow} C_\ast(\tilde
  k_q;X,X\setminus\eta_e)\,\cong\,C_{\ast-r}(k;P)\,.
\end{equation}

Let us now consider the following chain level definition of the
canonical isomorphism $H_\ast(f;X)\cong H_\ast(\tilde k_q;X)$, known
in Floer theory as the continuation isomorphism. Namely, given generic
Riemannian metrics both  for the negative gradient flow of $f$ and
of $\tilde k_q$, we have the chain complex morphism, which is a chain
homotopy equivalence,
$$
  \Phi\colon C_\ast(f;X)\,\to\,C_\ast(\tilde k_q;X),\quad
  m\mapsto \sum_{i(m')=i(m)}n(m,m')m',
$$
with
$$
  n(m,m')\,=\,\#_{\mathsf{alg}}\big(W^u(m;f)\pitchfork W^s(m';\tilde
  k_q)\big)\,.
$$
Composing $\Phi$ with the chain morphism from (\ref{eq collapse}), we
see that we obtain up to chain homotopy equivalence exactly the chain morphism $e_\bullet\colon
C_\ast(f;X)\to C_{\ast-r}(k;P)$.

We conclude the proof by comparing this construction in Morse homology
with the definition of the umkehr map in standard homology. We have
the commutative diagram
$$\begin{CD}
H_\ast(f;X)     @>\cong >>   H_\ast(\tilde k_q;X)  @>j_\ast>>
H_\ast(\tilde k_q;X,X\setminus\eta_e)   @>\cong>>  H_{\ast-r}(k;P)\\
@V\cong VV && @V\cong VV   @V\cong VV\\
H_\ast(X) & @>j_\ast >> & H_\ast (D(\nu_e),S(\nu_e)) @>{\cdot\cap
  u_{\nu_e}}>> H_{\ast-r}(P)\\
@V{\id}VV && @V\cong VV  @V{\id}VV\\
H_\ast(X) & @>{(\tau_e)_\ast}>> & H_\ast(P^{\nu_e}) @>\cong >>
H_{\ast-r}(P)
\end{CD}$$
where the upper row gives $e_!$ and the lower row is by definition
the umkehr map.
\end{proof}
This Morse-theoretic chain-level construction of the umkehr map provides Proposition 8 as an application of Proposition \ref{prop umkehr} to the proper embedding
$$
  r_{in}\colon Map\big(r(\G),M\big) \hookrightarrow LM^p\,.
  $$
  
      \vfill \eject

 \end{document}